\documentclass[11pt]{amsart}

\usepackage{amsthm, amsfonts, amssymb, amscd}
\usepackage[pagebackref,colorlinks]{hyperref}
\usepackage{tikz-cd}
\usepackage{geometry}
\usepackage{marginnote}

\theoremstyle{definition}
\newtheorem{ntn}{Notation}[section]

\theoremstyle{plain}
\newtheorem{lem}[ntn]{Lemma}
\newtheorem{prp}[ntn]{Proposition}
\newtheorem{thm}[ntn]{Theorem}
\newtheorem{cor}[ntn]{Corollary}

\theoremstyle{definition}
\newtheorem{question}[ntn]{Question}

\newtheorem{rem}[ntn]{Remark}
\newtheorem{exa}[ntn]{Example}

\numberwithin{equation}{section}

\newcommand{\N}{\mathbb{N}}
\newcommand{\z}{\mathbb{Z}}
\newcommand{\q}{\mathbb{Q}}

\newcommand{\C}{\mathbb{C}}

\newcommand{\Aa}{\mathcal{A}}
\newcommand{\BB}{\mathcal{B}}
\newcommand{\GG}{\mathcal{G}}

\newcommand{\WW}{\mathcal{W}}
\newcommand{\II}{\mathcal{I}}
\newcommand{\RP}{\mathcal{RP}}

\newcommand{\RB}{\mathcal{RB}}
\newcommand{\OO}{\mathcal{O}}
\newcommand{\KK}{\mathcal{K}}

\newcommand{\St}{\mathrm{{St}}}

\newcommand{\mmm}{\mathfrak{m}}
\newcommand{\ii}{\mathfrak{i}}

\renewcommand{\aa}{{A^\times}}

\newcommand{\tors}{{{\rm Tor}_1^{\z}}}

\newcommand{\mt}{\mapsto}
\newcommand{\lan}{\langle}
\newcommand{\ran}{\rangle}
\newcommand{\se}{\subseteq}
\newcommand{\arr}{\rightarrow}
\newcommand{\larr}{\longrightarrow}
\newcommand{\harr}{\hookrightarrow}
\newcommand{\two}{\twoheadrightarrow}
\newcommand{\Lan}{\langle\! \langle}
\newcommand{\Ran}{\rangle \!\rangle}

\newcommand{\stabe}{{\rm Stab}}
\newcommand{\Bb}{\mathrm{B}}
\newcommand{\Tt}{\mathrm{T}}
\newcommand{\Nn}{\mathrm{N}}
\newcommand{\GL}{\mathit{{\rm GL}}}

\newcommand{\Ee}{\mathit{{\rm E}}}
\newcommand{\PEe}{\mathit{{\rm PE}}}
\newcommand{\SL}{\mathit{{\rm SL}}}
\newcommand{\PSL}{\mathit{{\rm PSL}}}

\newcommand{\GW}{{\rm GW}}

\newcommand{\Ind}{{\rm Ind}}

\newcommand{\ab}{{\rm ab}}

\renewcommand{\char}{{\rm char}}

\newcommand{\coker}{{\rm coker}}
\newcommand{\im}{{\rm im}}

\newcommand{\inc}{{\rm inc}}
\newcommand{\id}{{\rm id}}
\newcommand{\GE}{{\rm GE}}

\newcommand{\CC}{\mathrm{C}}
\newcommand{\K}{\mathrm{K}}

\newcommand {\mtx}[2]
{\left(\!\!\!
\begin{array}{cc}
#1   \\
#2 
\end{array}
\!\!\!\right)}

\newcommand {\mtxx}[4]
{\left(\!
\begin{array}{cc}
\!\!#1 & \!\!#2   \\
\!\!#3 & \!\!#4
\end{array}\!\!\!
\right)
}

\newtheoremstyle{athm}
{}
{}
{\itshape}
{}
{\scshape}
{}
{.5em}
{\thmnote{#3}}
\theoremstyle{athm}

\begin{document}
\title[Homology groups of the elementary group of degree two]
{The low dimensional homology groups of the elementary group of degree two}
\author{Behrooz Mirzaii}
\author{Elvis Torres P\'erez}

\address{\sf 
Instituto de Ci\^encias Matem\'aticas e de Computa\c{c}\~ao (ICMC), Universidade de S\~ao Paulo, 
S\~ao Carlos, Brazil}
\email{bmirzaii@icmc.usp.br}
\address
{\sf Department of Sciences - Mathematics Section, Pontifical Catholic University of Peru (PUCP), Lima, Peru
%Faculty of Sciences, National University of Engineering (UNI), Lima, Peru
}
\email{etorresp@pucp.edu.pe
%elvis.torres.p@uni.pe
}

\begin{abstract}
In this article, we study the first, second, and third homology groups of the elementary group 
$\Ee_2(A)$, where $A$ is a commutative ring. In particular, we establish a refined Bloch--Wigner 
type exact sequence over a semilocal ring (subject to mild restrictions on its residue fields) 
in which either $-1 \in (\aa)^2$ or $|\aa/(\aa)^2| \leq 4$.

\end{abstract}
\maketitle

For a commutative ring $A$ with unit, let $\Ee_2(A)$ denote the subgroup of the general linear group 
$\GL_2(A)$ generated by elementary matrices. Clearly, $\Ee_2(A) \se \SL_2(A)$. We say that $A$ is 
a $\GE_2(A)$-ring if $\Ee_2(A) = \SL_2(A)$. 

In two previous papers \cite{B-E--2023} and \cite{B-E-2023} 
%and \cite{B-E2024} 
we studied the low-dimensional homology groups of the elementary group $\Ee_2(A)$. In this paper, 
we further investigate these homology groups.

The first main result of this paper concerns the first homology group of $\Ee_2(A)$. We show that 
for any commutative ring $A$, there is an exact sequence
\[
H_2(\Ee_2(A),\z)\arr H_1(Y_\bullet(A^2))_{\Ee_2(A)} \arr A/M \arr H_1(\Ee_2(A), \z) \arr 0,
\]
where $Y_\bullet(A^2)$ is a complex of $\GE_2(A)$-modules obtained from certain unimodular vectors 
in $A^2$ (see Theorem~\ref{H1}). Here $M$ is the additive subgroup of $A$ generated by $x(a^2 - 1)$ 
and $3(b + 1)(c + 1)$, where $x \in A$ and $a, b, c \in \aa$. The rightmost map is induced by the 
homomorphism $A \arr \Ee_2(A)$ sending $x$ to the elementary matrix $E_{12}(x)$. The group $A/M$ and 
its connection with the first homology of $\Ee_2(A)$ were first studied by Cohn in 
\cite{cohn1966} and \cite{cohn1968}.

Our second main result concerns the second homology group $H_2(\Ee_2(A), \z)$. For simplicity, let 
$A$ be a semilocal ring, where all its residue fields have more than three elements. Then it is not 
difficult to show that this homology group fits into the exact 
sequence
\[
H_2(\Bb(A), \z) \arr H_2(\Ee_2(A), \z) \arr I^2(A) \arr 0,
\]
where $\Bb(A)$ is the subgroup of $\Ee_2(A)$ consisting of all matrices whose lower left-hand 
entry is equal to zero. Here $I(A)$ denotes the fundamental ideal of $A$. More precisely, let 
$\GG_A$ be the square class group of $A$ and $\II_A$ the augmentation ideal of $\z[\GG_A]$. 
Let $\lan a \ran$ denotes the element of $\GG_A$ represented by $a \in \aa$. Then we define
\[
I(A) := \II_A / \lan \Lan a \Ran \Lan 1 - a \Ran : a, 1 - a \in \aa \ran,
\]
where $\Lan a \Ran := \lan a \ran - 1\in \z]\GG_A]$. Observe that $I^2(A) := (I(A))^2$.
Our second result concerns the map $H_2(\Ee_2(A), \z) \arr I^2(A)$ appearing in the above exact 
sequence. When $2 \in \aa$, we exhibit a cycle whose homology class in $H_2(\Ee_2(A), \z)$ maps 
to the element $\Lan a\Ran \Lan b\Ran$, for $a, b \in \aa$ (see Proposition~\ref{cycle} for the 
precise statement).

The third main result of this paper is a refined Bloch–Wigner exact sequence over a semilocal domain, 
which describes the third homology of $\Ee_2(A)$. More precisely, let $A$ be a semilocal ring such 
that for any maximal ideal $\mmm \in {\rm Specm}(A)$, the residue field $A/\mmm$ is either infinite, 
or if $|A/\mmm| = p^{d}$, then $(p - 1)d > 6$. If $-1 \in \aa^2$ or $|\aa/(\aa)^2| \leq 4$, then 
we obtain the refined Bloch–Wigner exact sequence
\[
0 \arr \tors(\mu(A),\mu(A))^\sim \arr 
\frac{H_3(\Ee_2(A),\z)}{\mu_2(A)\otimes_\z I^2(A)} 
\arr \RB(A) \arr 0,
\]
where $\RB(A)$ is the refined Bloch group of $A$ (see Theorem \ref{SL-PSL--1}). Here 
the map $\mu_2(A)\otimes_\z I^2(A) \arr H_3(\Ee_2(A),\z)$ is obtained from the product map 
\[
\mu_2(A)\otimes_\z H_2(\Ee_2(A),\z) \arr H_3(\Ee_2(A),\z),
\]
which is induced by the map $\mu_2(A)\times \Ee_2(A) \arr \Ee_2(A)$, $(-1, X)\mapsto -X$. This 
generalizes the main result of \cite{B-E-2023}. In \cite[Question~0.1]{B-E---2024}, 
we raised the question of whether, over local rings, the product map
\[
\mu_2(A)\otimes_\z H_2(\Ee_2(A),\z) \arr H_3(\Ee_2(A),\z)
\]
is trivial, provided that $|A/\mmm_A| \neq 2$.

Our main tool for studying the homology groups of $\Ee_2(A)$ is a first-quadrant
spectral sequence derived from the complex $Y_\bullet(A^2)$. All our main results are obtained 
from a detailed analysis of this spectral sequence. 

We now outline the organization of the present paper. In Section~\ref{s1} we give a brief overview
of the groups $\Ee_2(A) \subseteq \SL_2(A)$, $\GE_2(A) \subseteq \GL_2(A)$, the rank-one Steinberg 
group $\St(2,A)$, and the rank-one $\K_2$-group $\K_2(2,A)$. In Section~\ref{s2}, we introduce and 
study an $\SL_2(A)$-complex $L_\bullet(A^2)$ and its $\Ee_2(A)$-subcomplex $Y_\bullet(A^2)$. 
In particular, we study the groups $H_0$ and $H_1$ of these complexes. In Section~\ref{s3}, we 
introduce and analyze our main spectral sequence. In Section~\ref{s4}, we study the first homology 
group of $\Ee_2(A)$ and prove the first exact sequence discussed in this introduction. In 
Section~\ref{s5}, we study the natural map $H_2(\Ee_2(A), \z) \arr I^2(A)$. In Section~\ref{s6}, 
we establish the refined Bloch–Wigner exact sequence, which describes the third homology of 
$\Ee_2(A)$ over semilocal rings.
~\\
~\\
{\bf Notations.}
Throughout this paper all rings are commutative, except possibly group rings, and possess a unit 
element $1$. For a commutative ring $A$, let $\GG_A := \aa/(\aa)^2$ and define
\[
\WW_A := \{\, a \in A : a(1 - a) \in \aa \,\}.
\]
The element of $\GG_A$ represented by $a \in \aa$ is denoted by $\lan a \ran$, and we set 
$\Lan a\Ran := \lan a\ran - 1 \in \z[\GG_A]$. For $n \in \N$, let $\mu_n(A)$ denote the group of 
$n$th roots of unity in $A$. Moreover, define
\[
\mu(A) := \bigcup_{n \in \N} \mu_n(A).
\]
If $\BB \arr \Aa$ is a homomorphism of abelian groups, we write $\Aa / \BB$ for 
$\coker(\BB \arr \Aa)$.

%%%%%%%%%%%%%%%%%%%%%%%%%%%%%%%%%%%%%%%%%%%%%%%%%%%%%%%%%%%%%%%%%%%%%%%%%%
\section{Elementary and special linear groups of degree two}\label{s1}
%%%%%%%%%%%%%%%%%%%%%%%%%%%%%%%%%%%%%%%%%%%%%%%%%%%%%%%%%%%%%%%%%%%%%%%%%%

Let $A$ be a commutative ring. Let $\Ee_2(A)$ be the subgroup of $\GL_2(A)$ generated by the 
elementary matrices
$E_{12}(a):=\begin{pmatrix}
1 & a\\
0 & 1
\end{pmatrix}$ 
and 
$E_{21}(a):=\begin{pmatrix}
1 & 0\\
a & 1
\end{pmatrix}$, $a\in A$. Clearly 
\[
\Ee_2(A)\se \SL_2(A).
\]
For most commutative rings, $\Ee_2(A)$ is a proper subgroup of $\SL_2(A)$ and most of times it 
is even not a normal subgroup. We say that $A$ is a $\GE_2$-{\it ring} if 
\[
\Ee_2(A)=\SL_2(A).
\]

Let $\GE_2(A)$ be the subgroup of $\GL_2(A)$ generated by diagonal and elementary matrices. 
It is easy to verify that $\Ee_2(A)$ is normal in $\GE_2(A)$ 
%\cite[Proposition 2.1]{cohn1966} 
and the center of $\GE_2(A)$ is $\aa I_2$. Note that $A$ is a $\GE_2$-ring if and only if 
\[
\GE_2(A)=\GL_2(A).
\]

Semilocal rings and euclidean domains are examples of 
$\GE_2$-rings \cite[p. 245]{silvester1982}, \cite[\S 2]{cohn1966}. If $F$ is a field, then 
$F[X_1,\dots,X_d]$, $d>1$, is not a $\GE_2$-ring \cite[page 26]{cohn1966}.

A {\it quadratic ring} is a subring of $\C$ of the form $\z[\omega]:=\z \oplus \z\omega$, where 
$\omega$ is the root of an irreducible quadratic polynomial of the form $X^2\pm d$ or $X^2+X\pm d$ with 
$d$ a positive integer. The quadratic rings are split by the choice of sign into real quadratic rings
$\z[\sqrt{d}]$,  $\z[\frac{1}{2}(1+\sqrt{1+4d})]$ and imaginary quadratic rings $\z[\sqrt{-d}]$,
$\z[\frac{1}{2}(1+\sqrt{1-4d})]$. 
%Thus a quadratic ring either is of the form $\z[\sqrt{-d}]$ or of the form $\z[\frac{1}{2}(1+\sqrt{1-4d})]$.

\begin{thm}\label{non-normal}
Let $A=\z[\sqrt{-d}]$ or $A=\z[\frac{1}{2}(1+\sqrt{1-4d})]$, where $d$ is a positive integer.
\par {\rm (i)} If $d\geq 4$, then $\Ee_2(A)$ is a non-normal, infinite-index subgroup of $\SL_2(A)$.
\par {\rm (ii)} If $1\leq d < 4$, then $\Ee_2(A)=\SL_2(A)$, i.e. $A$ is a $\GE_2$-ring.
\end{thm}
\begin{proof}
See \cite[Theorem 1.5]{nica2011}.
\end{proof}

If $1 \leq d<4$, then $A$ is one of the rings $\z[\sqrt{-1}]$, $\z[\sqrt{-2}]$, $\z[\sqrt{-3}]$, 
$\z[\frac{1}{2}(1+\sqrt{-3})]$, $\z[\frac{1}{2}(1+\sqrt{-7})]$, $\z[\frac{1}{2}(1+\sqrt{-11})]$. 
Thus these rings are $\GE_2$-rings.

Let $K$ be an algebraic number field and let $\OO_K$ be its ring of algebraic integers. 
An {\it order} of $K$ is a subring $\OO$ of $\OO_K$ which is a free $\z$-module of rank $n=[K:\q]$. 
Thus any quadratic ring is an order. 
%Note that $\z[2i]:=\z\oplus 2\z(\sqrt{-1})$ is an order, which is a quadratic ring.

\begin{thm}[Liehl \cite{liehl1981}]\label{order}
Let A be an order in a number field which is not an imaginary quadratic ring. Then $\Ee_2(A)$
is a normal, finite-index subgroup in $\SL_2(A)$. Moreover, if the number field has a real 
embedding, then $\Ee_2(A)=\SL_2(A)$.  
\end{thm}
\begin{proof}
See \cite[Theorem 1.6]{nica2011}.
\end{proof}

The {\it rank one Steinberg group} $\St(2,A)$ of a commutative ring $A$ is the group with generators $x_{12}(r)$ 
and $x_{21}(s)$, $r,s\in A$, subject to the {\it Steinberg relations}
\begin{itemize}
\item[($\alpha$)] $x_{ij}(r)x_{ij}(s)=x_{ij}(r+s)$ for any $r,s \in A$,
\item[($\beta$)] $w_{ij}(u)x_{ji}(r)w_{ij}(u)^{-1}=x_{ij}(-u^2r)$, for any $u\in \aa$ and $r\in A$,\\
where $w_{ij}(u):=x_{ij}(u)x_{ji}(-u^{-1})x_{ij}(u)$.
\end{itemize}
The natural map $\St(2,A) \arr \Ee_2(A)$, given by $x_{ij}(r) \mapsto E_{ij}(r)$, is a well defined homomorphism of groups. 
The kernel of this map is called the {\it rank one $K_2$-group of} $A$ and is denoted by $\K_2(2,A)$:
\[
\K_2(2,A):=\ker (\St(2,A) \arr \Ee_2(A)).
\]
For $u\in \aa$, let 
\[
h_{ij}(u):=w_{ij}(u)w_{ij}(-1).
\]
An element of the form 
\[
\{v,u\}:=\{v,u\}_{12}=h_{12}(uv)h_{12}^{-1}(u)h_{12}(v)^{-1}
\]
is called a {\it Steinberg symbol}. This is an element of $\K_2(2,A)$ and is a central element of 
$\St(2,A)$ (see \cite[\S9]{D-S1973}). 

Let $\CC(2,A)$ be the subgroup of $\K_2(2,A)$ generated by the Steinberg symbols. A ring 
$A$ is called {\it universal for $\GE_2$} if $\K_2(2, A)$ is generated by Steinberg symbols 
\cite[p.~251]{D-S1973}. This definition coincides with the original definition of a ring 
universal for $\GE_2$ given by Cohn in \cite[page 8]{cohn1966} (for a proof of this fact see 
\cite[App. A]{h2022}).

%%%%%%%%%%%%%%%%%%%%%%%%%%%%%%%%%%%%%%%%%%%%%%%%%%%%%%%%%%%%%%%%%%%%%%%%%%
\section{The complex of unimodular vectors}\label{s2}
%%%%%%%%%%%%%%%%%%%%%%%%%%%%%%%%%%%%%%%%%%%%%%%%%%%%%%%%%%%%%%%%%%%%%%%%%%

A vector ${\pmb u}={\mtx {u_1} {u_2}}\in A^2$ is called {\it unimodular} if there exists a vector 
${\pmb v}={\mtx {v_1} {v_2}}\in A^2$ such that $({\pmb u},{\pmb v}):={\mtxx {u_1} {v_1} {u_2} {v_2}}\in\GL_2(A)$ 
and it is called $\GE_2$-{\it unimodular} if $({\pmb u}, {\pmb v})\in \GE_2(A)$.

For any non-negative integer $n$, let $L_n(A^2)$ be the free abelian group generated by 
the set of all $(n+1)$-tuples $(\lan{\pmb v_0}\ran, \dots, \lan{\pmb v_n}\ran)$, where every vector
${\pmb v_i}\in A^2$ is unimodular and for any $i\neq j$, $({\pmb v_i},{\pmb v_j})\in\GL_2(A)$. 
Let $Y_n(A^2)$ be the free abelian subgroup of $L_n(A^2)$ generated by the set of all 
$(n+1)$-tuples $(\lan{\pmb v_0}\ran, \dots, \lan{\pmb v_n}\ran)$, where every vector
${\pmb v_i}\in A^2$ is $\GE_2$-unimodular. It is easy to see if ${\pmb u}$ is
$\GE_2$-unimodular and $({\pmb u},{\pmb v})\in\GL_2(A)$, then $({\pmb u},{\pmb v})\in\GE_2(A)$
and hence ${\pmb v}$ is $\GE_2$-unimodular.

We consider $L_n(A^2)$ (resp. $Y_n(A^2)$) as a left $\SL_2(A)$-module (resp. left $\Ee_2(A)$-module) in 
a natural way. If necessary, we convert this action to a right action by the definition
$m.g:=g^{-1}m$. Let us define the $n$-th differential operator
\[
\partial_n^L : L_n(A^2) \arr L_{n-1}(A^2), \ \ n\ge 1,
\]
as an alternating sum of face operators which throws away the $i$-th component of generators. Hence we 
have the complex $L_\bullet(A^2)$ with the  $\Ee_2(A)$-subcomplex $Y_\bullet(A^2)$.

For a subgroup $H$ of a group $G$ and any $H$-module $M$, let $\Ind_H^GM:=\z[G]\otimes_H M$ be the 
induced module (from $H$ to $G$). This extension of scalars is called {\it induction} from $H$ to $G$ 
\cite[\S5, Chap. III]{brown1994}.

\begin{prp}\label{ind1}
For any commutative ring $A$, 
%\[
$L_\bullet(A^2)\simeq \Ind_{\Ee_2(A)}^{\SL_2(A)}Y_\bullet(A^2)$.
%\]
\end{prp}
\begin{proof}
Clearly the map
\[
\varphi_\bullet: \z[\SL_2(A)] \otimes_{\Ee_2(A)} Y_\bullet(A^2) \arr L_\bullet(A^2)
\]
defined by
\[
g\otimes (\lan\pmb{v_0}\ran, \lan\pmb{v_1}\ran, \dots, \lan\pmb{v_n}\ran) \mapsto 
(\lan g\pmb{v_0}\ran, \lan g\pmb{v_1}\ran, \dots, \lan g\pmb{v_n}\ran),
\]
is a well-defined morphism of $\SL_2(A)$-complexes. This is an 
isomorphism with the inverse morphism
\[
\psi_\bullet: L_\bullet(A^2) \arr \z[\SL_2(A)] \otimes_{\Ee_2(A)} Y_\bullet(A^2)
\]
defined by
\[
(\lan\pmb{v_0}\ran, \lan\pmb{v_1}\ran, \dots, \lan\pmb{v_n}\ran) \mapsto g \otimes 
(\lan\pmb{e_1}\ran, \lan g^{-1}\pmb{v_1}\ran, \dots, \lan g^{-1}\pmb{v_n}\ran),
\]
where $g\pmb{e_1}=\pmb{v_0}$ for some $g \in \SL_2(A)$. 
\end{proof}

\begin{cor}\label{ind2}
For any non-negative integer $n$, 
%\[
$H_n(L_\bullet(A^2))\simeq \Ind_{\Ee_2(A)}^{\SL_2(A)} H_n(Y_\bullet(A^2))$.
%\]
\end{cor}
\begin{proof}
Since the functor $\Ind_{\Ee_2(A)}^{\SL_2(A)}$ is exact on the category of $\z[\Ee_2(A)]$-modules (since
$\z[\SL_2(A)]$ is a free $\z[\Ee_2(A)]$-module), the claim follows from the previous lemma.
\end{proof}

Let $\partial_{-1}=\varepsilon: L_0(A^2) \arr \z$ be the natural augmentation map, which is given by 
$\sum_i n_i(\lan {\pmb v_{0,i}}\ran) \mt \sum_i n_i$. Then we have the complexes
$L_\bullet(A^2)\overset{\varepsilon}{\arr} \z$ and $Y_\bullet(A^2)\overset{\varepsilon}{\arr} \z$.

\begin{thm}[Hutchinson]\label{GE2C}
For any commutative ring $A$, $H_0(Y_\bullet(A^2))\overset{\bar{\varepsilon}}{\simeq} \z$, where $\bar{\varepsilon}$
is induced by $\varepsilon$. In other words, 
the complex $Y_1(A^2) \overset{\partial_1^Y}{\larr} Y_0(A^2) \overset{\varepsilon}{\arr} \z \arr 0$ 
always is exact. Moreover,
\[
H_0(L_\bullet(A^2))\simeq \Ind_{\Ee_2(A)}^{\SL_2(A)} \z\simeq \z[\SL_2(A)/\Ee_2(A)],
\]
where the right side is the coset module.
\end{thm} 
\begin{proof}
For the first isomorphism see the proof of \cite[Theorem 3.3]{h2022}. The isomorphism 
$H_0(L_\bullet(A^2))\simeq \Ind_{\Ee_2(A)}^{\SL_2(A)} \z$ follows from 
the first isomorphism and Corollary \ref{ind2}. The other isomorphism follows from the fact that
$\z[\SL_2(A)]$ is a free $\z[\Ee_2(A)]$-module with  a basis consisting of one representative 
for each element of the set $\SL_2(A)/\Ee_2(A)$.

\end{proof}

\begin{exa}
(i) Let $A=\z[\sqrt{-d}]$ or $A=\z[\frac{1}{2}(1+\sqrt{1-4d})]$, where $d\geq 4$. Then by Theorem \ref{GE2C}
and Theorem \ref{non-normal}, $H_0(L_\bullet(A^2))$ is a free abelian group of infinite and countable rank.
\par(ii) If $A$ is an order in a number field which is not an imaginary quadratic ring, then $\Ee_2(A)$
is a normal and finite-index subgroup of $\SL_2(A)$ (see Theorem~\ref{order}). Hence $H_0(L_\bullet(A^2))$
is a free abelian group of finite rank.
%\par (iii) Let $A$ be a ring. Then by Theorem \ref{GE2C},
%\[
%H_0(L_\bullet(A^2))\simeq \Ind_{\Ee_2(A)}^{\SL_2(A)} \z\simeq \z[\SL_2(A)/\Ee_2(A)],
%\]
%where the right side is the coset module.
\end{exa}

\begin{thm}[Hutchinson]\label{H1(Y)}
For any commutative ring $A$, we have the isomorphisms
\[
H_1(Y_\bullet(A^2))\simeq \bigg(\displaystyle{\frac{\K_2(2, A)}{\CC(2, A)}}\bigg)^\ab, \ \ \ \ \
H_1(L_\bullet(A^2))\simeq \Ind_{\Ee_2(A)}^{\SL_2(A)}\bigg(\displaystyle{\frac{\K_2(2, A)}{\CC(2, A)}}\bigg)^\ab.
\]
In particular if $A$ is universal for $\GE_2$, then $Y_\bullet(A^2) \overset{\epsilon}{\arr} \z$ and 
$L_\bullet(A^2) \overset{\epsilon}{\arr} \z$ are exact in  dimension $1$.
\end{thm}
\begin{proof}
The first isomorphism follows from \cite[Theorem 7.2]{h2022}. The second isomorphism follows from the first 
isomorphism and Proposition \ref{ind2}.
\end{proof}

It follows from Lemma \ref{ind1} and Shapiro’s Lemma that the inclusion 
$Y_\bullet(A^2)\arr L_\bullet(A^2)$ induces isomorphisms of the homology groups
\[
H_q(\Ee_2(A), Y_p(A^2)) \simeq H_q(\SL_2(A), L_p(A^2))
\]
for all $ p,q$. For any $n\geq 0$, let
\[
Z_n^{Y}(A^2):=\ker(\partial_n^Y), \ \ \ \ \ Z_n^{L}(A^2):=\ker(\partial_n^L).
\]

\begin{lem}\label{ind3}
For all $p,q\geq 0$, we have the isomorphism 
\[
H_q(\Ee_2(A),  Z_p^{Y}(A^2)) \simeq H_q(\SL_2(A), Z_p^{L}(A^2)).
\]
\end{lem}
\begin{proof}
The functor $\Ind_{\Ee_2(A)}^{\SL_2(A)}$ is exact on the 
category of $\z[\Ee_2(A)]$-modules. Thus we have the isomorphism
\[
Z_n^{L}(A^2) \simeq \Ind_{\Ee_2(A)}^{\SL_2(A)} Z_n^{Y}(A^2).
\]
Now the claim follows by applying Shapiro's Lemma to this isomorphism.
\end{proof}

Let $G$ be a group and let $C_\bullet$ be a complex of left $G$-modules: 
$C_\bullet: \ \ \cdots \arr C_1 \arr C_0\arr 0$.
The $n$-th homology of $G$ with coefficients in $C_\bullet$, denoted by 
\[
H_n(G, C_\bullet),
\]
is defined as the $n$-th homology of the total complex of the double complex 
$F_\bullet \otimes_{G} C_\bullet$, where $F_\bullet\arr \z$ is a projective 
resolution of $\z$ over $G$ \cite[\S3, Chap. VII]{brown1994}.

By applying Shapiro's Lemma to the isomorphism of Proposition \ref{ind1},
for any non-negative integer $n$, we obtain the isomorphism
\[
H_n(\Ee_2(A),Y_\bullet(A^2)) \simeq H_n(\SL_2(A),L_\bullet(A^2)).
\]
Let $\z_\bullet$ be the complex with terms $\z_0=\z$ and $\z_n=0$ for $n\neq 0$. Then $\varepsilon$ induces
the natural commutative diagram of complexes
\[
\begin{tikzcd}
Y_\bullet(A^2) \ar[rr, hook] \ar[dr, "\varepsilon"] && L_\bullet(A^2) \ar[dl, "\varepsilon"]\\
 & \z_\bullet&
\end{tikzcd}
\]
From this diagram we obtain the commutative diagram
\[
\begin{tikzcd}
H_n(\Ee_2(A), Y_\bullet(A^2)) \ar[r, "\simeq "] \ar[d, "\varepsilon"] & 
H_n(\SL_2(A), L_\bullet(A^2)) \ar[d, "\varepsilon"]\\
H_n(\Ee_2(A), \z) \ar[r] & H_n(\SL_2(A), \z).
\end{tikzcd}
\]
In Theorem \ref{H1} below we give a satisfactory description of $H_1(\Ee_2(A),\z)$ over any commutative ring, while
the study of $H_1(\SL_2(A),\z)$ seems much harder (see for example \cite{swan1971} and \cite{gmv1994}).

%%%%%%%%%%%%%%%%%%%%%%%%%%%%%%%%%%%%%%%%%%%%%%%%%%%%%%%%%%%%%%%%%%%%%%%%%%
\section{The main spectral sequence and the homology of \texorpdfstring{$\Ee_2$}{Lg}}\label{s3}
%%%%%%%%%%%%%%%%%%%%%%%%%%%%%%%%%%%%%%%%%%%%%%%%%%%%%%%%%%%%%%%%%%%%%%%%%%

Let $G$ be a group and let $C_\bullet$ be a complex of left $G$-modules. Let $B_\bullet(G)\arr \z$ 
be the bar resolution of $G$. From the double complex 
\[
B_\bullet(G) \otimes_{G} C_\bullet
\]
we obtain two spectral sequences
\[
\mathsf{E}_{p,q}^2(G)=H_p(G, H_q(C_\bullet))\Rightarrow H_{p+q}(G, C_\bullet),
\]
and 
\[
\mathrm{E}_{p, q}^1(G)=H_q(G, C_p) \Rightarrow H_{p+q}(G, C_\bullet),
\]
(see \cite[\S 5, Chap. VII]{brown1994}).

\begin{lem}\label{H0}
Let the complex $C_\bullet$ be exact for $1\leq  i \leq n$ and $M:=H_0(C_\bullet)$. Then 
$H_i(G, C_\bullet)\simeq H_i(G, M)$ for $0 \leq  i \leq n$. 
\end{lem}
\begin{proof}
Apply \cite[Chap.~VII, Proposition~5.2]{brown1994} to the morphism of complexes 
$C_\bullet \arr M_\bullet$, where $M_\bullet$ is the complex consisting of the single module 
$M$ concentrated in degree~$0$.
%This follows from an easy analysis of the spectral sequence $\mathsf{E}_{p,q}^2(G)$.
%\[
%\mathsf{E}_{p,q}^2(G)=H_p(G, H_q(L_\bullet))\Rightarrow H_{p+q}(G, L_\bullet).
%\]
\end{proof}

Thus if $C_\bullet$ is exact for every $i\geq 1$, then we have the spectral sequence
\begin{equation}\label{esp}
\mathrm{E}_{p, q}^1(G)=H_q(G, C_p) \Rightarrow H_{p+q}(G, H_0(C_\bullet)).
\end{equation}

For any commutative ring $A$, the sequence
\[
0 \arr Z_1^Y(A^2) \overset{\inc}{\arr} Y_1(A^2) \arr Y_0(A^2) \overset{\varepsilon}{\arr} \z \arr 0,
\]
is exact (Theorem \ref{GE2C}). 
%From the complex 
%\[
%0 \arr Z_1^Y(A^2) \overset{\inc}{\arr} Y_1(A^2) \overset{\partial_1^Y}{\arr} Y_0(A^2) \arr 0,
%\]
%we obtain the double complex
%\[ 
%0 \! \arr \! B_\bullet(\Ee_2(A)) \! \otimes_{\Ee_2(A)} \! Z_1^Y(A^2) \! \arr \! B_\bullet(\Ee_2(A))
%\!\otimes_{\Ee_2(A)} \! Y_1(A^2) \! \arr \! B_\bullet(\Ee_2(A)) \! \otimes_{\Ee_2(A)} \! Y_0(A^2) \!\arr \! 0,
%\]
%where $B_\bullet(\Ee_2(A)) \arr \z$ is the bar resolution of $\z$ over $\Ee_2(A)$.
Thus by (\ref{esp}) we have the first quadrant spectral sequence
\[
E^1_{p.q}=\left\{\begin{array}{ll}
H_q(\Ee_2(A),Y_p(A^2)) & p=0,1\\
H_q(\Ee_2(A),Z_1^Y(A^2)) & p=2\\
0 & p>2
\end{array}
\right.
\Longrightarrow H_{p+q}(\Ee_2(A),\z).
\]
%(see \cite[§5, Chap. VII]{brown1994} or \cite[\S1]{mirzaii2017}).

The group $\Ee_2(A)$ acts transitively on the sets of generators of $Y_i(A^2)$ for $i=0,1$. Let
\[
{\pmb\infty}:=\lan {\pmb e_1}\ran, \ \ \  {\pmb 0}:=\lan {\pmb e_2}\ran , \ \ \  
{\pmb a}:=\lan {\pmb e_1}+ a{\pmb e_2}\ran, \ \ \ a\in \aa,
\]
where ${\pmb e_1}:={\mtx 1 0}$ and $ {\pmb e_2}:={\mtx 0 1}$.
We choose $({\pmb \infty})$ and $({\pmb \infty} ,{\pmb 0})$ as 
representatives of the orbit of the generators of $Y_0(A^2)$ and $Y_1(A^2)$,
respectively. Then 
\[
Y_0(A^2)\simeq \Ind _{\Bb(A)}^{\Ee_2(A)}\z, \ \ \ \ \ \ \ \ \ \ \ \ 
Y_1(A^2)\simeq \Ind _{\Tt(A)}^{\Ee_2(A)}\z,
\]
where 
\[
\Bb(A):=\stabe_{\Ee_2(A)}({\pmb \infty})=\Bigg\{\begin{pmatrix}
a & b\\
0 & a^{-1}
\end{pmatrix}:a\in \aa, b\in A\bigg\},
\]
\[
\Tt(A):=\stabe_{\Ee_2(A)}({\pmb \infty},{\pmb 0})=\Big\{D(a):a\in \aa\Big\},
\]
where 
$D(a):=\begin{pmatrix}
a & 0\\
0 & a^{-1}
\end{pmatrix}$ (see \cite[\S3, Chap. I]{brown1994}).
Note that $\Tt(A)\simeq \aa$ and $\Bb(A)\simeq A \rtimes \aa$. In our calculations
sometimes we identify $\Tt(A)$ with $\aa$. By Shapiro's lemma we have
\[
\mathrm{E}_{0,q}^1 \simeq H_q(\Bb(A),\z), \ \ \ \ \ \ \mathrm{E}_{1,q}^1 \simeq H_q(\Tt(A),\z).
\]
In particular, $\mathrm{E}_{0,0}^1\simeq \z\simeq \mathrm{E}_{1,0}^1$. Moreover, $d_{1, q}^1=H_q(\sigma) - H_q(\inc)$,
where $\sigma: \Tt(A) \arr \Bb(A)$ is given by $\sigma(D(a))= wD(a) w^{-1}=D(a^{-1})$ for 
$w={\mtxx 0 1 {-1} 0}$. This easily implies that  $d_{1,0}^1$ is trivial and 
\[
d_{1,1}^1: \Tt(A) \arr \Bb(A)^\ab=H_1(\Bb(A),\z), \ \ \ \ \ D(a)\mt D(a)^{-2}.
\]
In fact, note that 
\[
\Bb(A)^\ab\simeq \Tt(A)\oplus A_\aa\simeq \aa\oplus A_\aa,
\]
where 
\[
A_\aa=A/\lan a^2-1: a\in \aa \ran.
\]
Thus $d_{1,1}^1(D(a))=(D(a)^{-2}, 0)$.  This implies that 
\[
\mathrm{E}_{0,1}^2\simeq\GG_A\oplus A_\aa, \ \ \ \ker(d_{1,1}^1)\simeq\mu_2(A).
\]
Furthermore, it is not difficult to see that $d_{2,1}^1:H_1(\Ee_2(A),Z_1^Y(A^2)) \arr \mu_2(A)$
is surjective. In fact, for any $b\in \mu_2(A)$, we have
$d_{2,1}^1([b]\otimes \partial_2^Y({\pmb \infty}, {\pmb 0}, {\pmb a}))=b$. Therefore 
\[
\mathrm{E}_{1,1}^2=0.
\]

The orbits of the action of  $\Ee_2(A)$ on $Y_2(A)$ are represented by
\[
\lan a\ran[\ ]:=({\pmb\infty}, {\pmb 0},{\pmb a}), \ \ \ \ \lan a\ran\in \GG_A.
\]
Therefore $Y_2(A^2)_{\Ee_2(A)} \simeq \z[\GG_A]$. Moreover, the map
\[
\z[\GG_A]\simeq Y_2(A^2)_{\Ee_2(A)} \arr Y_1(A^2)_{\Ee_2(A)}\simeq \z
\]
coincides with the augmentation map. Now from the composite 
\[
\z[\GG_A]\simeq Y_2(A^2)_{\Ee_2(A)} \arr Z_1^Y(A^2)_{\Ee_2(A)} \arr Y_1(A^2)_{\Ee_2(A)}\simeq \z
\]
it follows that $d_{2,0}^1:H_0(\Ee_2(A),Z_1^Y(A^2)) \arr \z$ is surjective. Hence $\mathrm{E}_{1,0}^2=0$. Let 
\[
\GW(A):=H_0(\Ee_2(A), Z_1^Y(A^2)),  \ \ \ \ \epsilon:=d_{2,0}^1:\GW(A) \arr \z.
\]
We denote the kernel of $\epsilon$  by $I(A)$. 

Combining all these results, we obtain that the first page of the spectral sequence has the following form: \\[2mm]
\[
\begin{tikzcd}
H_2(\Bb(A),\z) & \ar[l,"0"'] \aa \wedge \aa & \ar[l,"d_{2,2}^1"'] H_2(\SL_2(A), Z_1^Y(A^2)) & 0 \\
\aa \oplus A_\aa & \ar[l,"d_{1,1}^1"'] \aa & \ar[l,"d_{2,1}^1"'] H_1(\SL_2(A), Z_1^Y(A^2)) & 0 \\
\z & \ar[l,"0"'] \z & \ar[l,"d_{2,0}^1"', two heads] \GW(A) & 0
\end{tikzcd}
\]
The second page of the spectral sequence is given by:
\[
\begin{tikzcd}
H_2(\Bb(A),\z) & \mathrm{E}_{1,2}^2 & \mathrm{E}_{2,2}^2 & 0 \\
\GG_A \oplus A_{\aa} & 0 & \ar[llu,"d^2_{2,1}"'] \mathrm{E}_{2,1}^2 & 0 \\
\z & 0 & \ar[llu,"d^2_{2,0}"'] I(A) & 0
\end{tikzcd}
\]

Let 
\[
\overline{\GW}(A):=\z[\GG_A]/\lan \Lan a\Ran \Lan 1-a\Ran : a\in \WW_A\ran. 
\]
The augmentation map $\z[\GG_A]\arr\z$ induces the natural map $\bar{\epsilon}:\overline{\GW}(A)\arr\z$.
The kernel of this map is denoted by $\bar{I}(A)$. Note that 
\[
\bar{I}(A):=\II_A/\lan \Lan a\Ran \Lan 1-a\Ran: a\in \WW_A \ran,
\]
where $\II_A$ is the augmentation ideal of $\z[\GG_A]$.

From the exact sequence $0 \arr B_1^Y(A^2) \arr Z_1^Y(A^2) \arr H_1(Y_\bullet(A^2)) \arr 0$,
we obtain the exact sequence 
\[
H_0(\Ee_2(A), B_1^Y(A^2)) \arr \GW(A) \arr H_1(Y_\bullet(A^2))_{\Ee_2(A)}\arr 0.
\]
From the sequence $Y_3(A^2)\arr Y_2(A^2) \two B_1^Y(A^2)$ we obtain the sequence
\[
H_0(\Ee_2(A), Y_3(A^2)) \arr H_0(\Ee_2(A), Y_2(A^2)) \two H_0(\Ee_2(A), B_1^Y(A^2)).
\]
We showed in above that $Y_2(A^2)_{\Ee_2(A)} \simeq \z[\GG_A]$.
The orbits of the action of $\Ee_2(A)$ on $Y_3(A)$ are  represented by
\[
\lan a\ran[x]:=({\pmb\infty}, {\pmb 0},{\pmb a}, \pmb{ax}),
\ \ \lan a\ran\in \GG_A, x \in \WW_A.
\]
Thus $Y_3(A^2)_{\Ee_2(A)} $ is the free $\z[\GG_A]$-module generated by 
the symbols $[x]$, $x\in \WW_A$. Moreover, the map
\[
Y_3(A^2)_{\Ee_2(A)} \arr Y_2(A^2)_{\Ee_2(A)}\simeq \z[\GG_A]
\]
is given by $[x]\mapsto \Lan x\Ran \Lan 1-x\Ran$. All these induce a surjective map
\[
\overline{\GW}(A) \two H_0(\Ee_2(A), B_1^Y(A^2)).
\]
Thus we have the exact sequence
\[
\overline{\GW}(A) \arr \GW(A) \arr H_1(Y_\bullet(A^2))_{\Ee_2(A)}\arr 0.
\]
Finally from this we obtain the exact sequence
\begin{align}\label{I-I}
\bar{I}(A) \arr I(A) \arr H_1(Y_\bullet(A^2))_{\Ee_2(A)}\arr 0.
\end{align}

%%%%%%%%%%%%%%%%%%%%%%%%%%%%%%%%%%%%%%%%%%%%%%%%%%%%%%%%%%%%%%%%%%%%%%%%%%
\section{The first homology of \texorpdfstring{$\Ee_2$}{Lg}}\label{s4}
%%%%%%%%%%%%%%%%%%%%%%%%%%%%%%%%%%%%%%%%%%%%%%%%%%%%%%%%%%%%%%%%%%%%%%%%%%

In this section, we show that the first homology group of $\Ee_2(A)$ fits into a natural 
exact sequence that almost completely describes its structure.

\begin{thm}\label{H1}
Let $A$ be a commutative ring and $M$ the additive subgroup of $A$ generated by $x(a^2-1)$ 
and $3(b+1)(c+1)$, where  $x\in A$ and $a,b,c \in \aa$. Then there is an exact sequence of 
$\GG_A$-modules
\[
H_2(\Ee_2(A),\z)\arr H_1(Y_\bullet(A^2))_{\Ee_2(A)} \arr A/M \arr H_1(\Ee_2(A), \z) \arr 0. 
\]
\end{thm}
\begin{proof}
By an easy analysis of the  above second page of the spectral sequence $E_{\bullet,\bullet}$, 
we obtain the exact sequence
\[
H_2(\Bb(A),\z)\arr  H_2(\Ee_2(A), \z) \arr I(A) \overset{d_{2,0}^2}{\larr} \GG_A\oplus A_\aa 
\arr H_1(\Ee_2(A),\z) \arr 0.
\]
Let's study the composite map $\bar{I}(A)\arr I(A) \overset{d_{2,0}^2}{\larr} \GG_A\oplus A_\aa$.
Take the element $\Lan a\Ran\in \bar{I}(A)$, $a\in \aa$. This is represented by 
\[
X_a:=[\ ]\otimes \Big((\pmb{\infty},\pmb{0},\pmb{a})-(\pmb{\infty},\pmb{0},\pmb{1})\Big)\in 
B_0(\Ee_2(A))\otimes_{\Ee_2(A)} Y_2(A^2).
\]
By passing to $I(A)$, the image of $\Lan a\Ran$ is represented by
\[
[\ ]\otimes \partial_2((\pmb{\infty},\pmb{0},\pmb{a})-(\pmb{\infty},\pmb{0},\pmb{1}))\in 
B_0(\Ee_2(A))\otimes_{\Ee_2(A)} Z_1(A^2),
\]
which we denote again by $X_a$. Thus $d^2_{2,0}(\Lan a \Ran)=d^2_{2,0}(\overline{X_a})$. 
To calculate the above composition we look at the diagram
\[
\begin{tikzcd}
B_1(\Ee_2(A))\otimes Y_0(A^2) & \ar[l,"\id\otimes\partial_1"'] 
B_1(\Ee_2(A))\otimes Y_1(A^2) \ar[d,"d_1\otimes\id"]\\
& B_0(\Ee_2(A))\otimes Y_1(A^2) & \ar[l,"\id\otimes\inc"'] B_0(\Ee_2(A))\otimes Z_1(A^2).
\end{tikzcd}
\]
We have
\begin{align*}
(\id_{B_0(\Ee_2(A))}\otimes\inc)(X_a)& =[\ ]\otimes\Big((\pmb{0},\pmb{a})-(\pmb{\infty},\pmb{a})
-(\pmb{0},\pmb{1})+(\pmb{\infty},\pmb{1})\Big)\\
&=[\ ]\otimes\Big(g_a(\pmb{\infty},\pmb{0})-h_a(\pmb{\infty},\pmb{0})
-g_1(\pmb{\infty},\pmb{0})+h_1(\pmb{\infty},\pmb{0})\Big)\\
&=(g^{-1}_a-h^{-1}_a-g^{-1}_1+h^{-1}_1)[\ ]\otimes(\pmb{\infty},\pmb{0}),
%\in B_0(\Ee_2(A))\otimes X_1(A^2)
\end{align*}
where $g_a=\mtxx{0}{1}{-1}{a}$, $h_a=\mtxx{1}{a^{-1}}{0}{1}$.
This element sits in the image of $d_1\otimes\id_{Y_1(A^2)}$. In fact
\[
(d_1\otimes\id_{Y_1(A^2)})\Big(([g^{-1}_a]-[h^{-1}_a]-[g^{-1}_1]+[h^{-1}_1])
\otimes(\pmb{\infty},\pmb{0})\Big)=
(g^{-1}_a-h^{-1}_a-g^{-1}_1+h^{-1}_1)[\ ]\otimes(\pmb{\infty},\pmb{0}).
\]
Finally $(\id_{B_1(\Ee_2(A))}\otimes\partial_1)\Big(([g^{-1}_a]-[h^{-1}_a]-[g^{-1}_1]+[h^{-1}_1])
\otimes(\pmb{\infty},\pmb{0})\Big)$ 
represents $d_{2,0}^ 2(\Lan a\Ran)$. This is the element
\[
\bigg(w\Big([g^{-1}_a]-[h^{-1}_a]-[g^{-1}_1]+[h^{-1}_1]\Big)-\Big([g^{-1}_a]
-[h^{-1}_a]-[g^{-1}_1]+[h^{-1}_1]\Big)\bigg)
\otimes(\pmb{\infty})\in B_1(\Ee_2(A))\otimes Y_0(A^2).
\] 
By adding the null element
$(d_2\otimes\id_{Y_0(A^2)})(([w|h^{-1}_a]-[w|g^{-1}_a]+[w|g^{-1}_1]-[w|h^{-1}_1])\otimes(\pmb{\infty}))$, 
we see that $d_{2,0}^ 2(\Lan a\Ran)$ is represented by
\[
([wg^{-1}_a]-[wh^{-1}_a]-[wg^{-1}_1]+[wh^{-1}_1]-[g^{-1}_a]+[h^{-1}_a]+[g^{-1}_1]-[h^{-1}_1])
\otimes(\pmb{\infty}).
\]
Now by adding the null element 
$(d_2\otimes\id_{Y_0(A^2)})(([D(a)^{-1}h^{-1} _{a^{-1}}|wh^{-1}_a]-[h^{-1} _1|wh^{-1}_1])\otimes(\pmb{\infty}))$ 
(note that $D(z)^{-1}h^{-1} _{z^{-1}}\in \Bb(A)=\stabe_{\Ee_2(A)}(\pmb{\infty})$ and 
$D(z)^{-1}h^{-1}_{a^{-1}}wh^{-1}_a=wg^{-1}_a$), we see that $d_{2,0}^ 2(\Lan a\Ran)$ is represented by
\[
([D(a)^{-1}h^{-1}_{a^{-1}}]-[h^{-1}_1]-[g^{-1}_a]+[h^{-1}_a]+[g^{-1}_1]-[h^{-1}_1])\otimes(\pmb{\infty}).
\]
Once more, by adding the null element 
\[
(d_2\otimes\id_{X_0(A^2)})\Big(([h_1|-w]-[h_{a^{-1}}|-w])\otimes(\pmb{\infty})\Big)
=([g^{-1}_a]-[h_{a^{-1}}]-[g^{-1}_1]+[h_1])\otimes(\pmb{\infty}),
\]
we see that $d_{2,0}^ 2(\Lan a\Ran)$ is represented by
\[
([D(a)^{-1}h^{-1}_{a^{-1}}]-[h^{-1}_1]+[h^{-1}_a]-[h^{-1}_1]-[h_{a^{-1}}]+[h_1])\otimes(\pmb{\infty}).
\]
Note that in $H_1(\Ee_2(A),Y_0(A^2))$, if $g,h\in \Bb(A)$, then 
\[
[gh]\otimes(\pmb{\infty})=([g]+[h])\otimes(\pmb{\infty}) \ \ 
\text{and} \ \ 
[h^{-1}]\otimes(\pmb{\infty})=-[h]\otimes(\pmb{\infty}).
\]
Thus $d_{2,0}^ 2(\Lan a\Ran)$ in $H_1(\Ee_2(A),Y_0(A^2))/\im(d^1_{1,1})$ is represented by
\[
[D(a)^{-1}h^{-2}_{a^{-1}}h^{-1}_ah^3_1]\otimes(\pmb{\infty})
=\left[\mtxx{a^{-1}}{0}{0}{a}\mtxx{1}{-2a-a^{-1}+3}{0}{1}\right]\otimes(\pmb{\infty}).
\]
Therefore in  $\GG_A\oplus A_{\aa}$ we have
\begin{align*}
d_{2,0}^ 2(\Lan a\Ran)&=(\lan a^{-1}\ran, \overline{-2a-a^{-1}+3})
=(\lan a\ran, \overline{3(1-a)})=3(\lan a \ran, \overline{1-a}).    
\end{align*}
Thus we get the exact sequence
\[
H_2(\Ee_2(A),\z) \arr I(A)/\bar{I}(A) \arr \frac{\GG_A \oplus A_\aa}
{\lan (\lan b\ran, 3(b-1))\ |b\in \aa \ran} \arr H_1(\Ee_2(A),\z) \arr 0.
\]
The map
\[
\frac{A}{M} \arr \frac{\GG_A \oplus A_\aa} {\lan (\lan b\ran, 3(b-1))\ |b\in \aa \ran}, 
\ \ \ \overline{a} \mt \overline{(\lan 1\ran, \bar{a})},
\]
is an isomorphism with the inverse map
\[
\frac{\GG_A \oplus A_\aa}
{\lan (\lan b\ran, 3(b-1))\ |b\in \aa \ran} \arr \frac{A}{M}, \ \ \ 
\overline{(\lan b\ran, \bar{a})} \mt \overline{a-3b+3}.
\]
Now (\ref{I-I}) completes the proof of the theorem.
\end{proof}

\begin{cor}
If $A$ is universal for $\GE_2$, then $H_1(\Ee_2(A),\z)\simeq A/M$, where $M$ is the additive 
subgroup of $A$ generated by $x(a^2-1)$ and $3(b+1)(c+1)$ with  $x\in A$ and $a,b,c\in\aa$.
\end{cor}
\begin{proof} 
By Theorem \ref{H1(Y)}, $H_1(Y_\bullet(A^2))=0$. Now the claim follows from the previous theorem.
\end{proof}

\begin{exa}\label{universal}
(i) Let $R$ be a local ring such that its residue field has at least four elements. 
Then there is $a\in R^\times$ such that $a^2-1\in R^\times$. Now for any $R$-algebra 
$A$ and for any $x\in A$, we have
\[
x=x(a^2-1)^{-1}(a^2-1)\in M.
\]
This shows that $A=M$ and therefore by the above theorem
\[
H_1(\Ee_2(A),\z)=0.
\]
\par (ii) Let $A$ be a local ring with maximal ideal $\mmm_A$. Then by the above corollary
%Theorem~\ref{H1(Y)} and Theorem~\ref{H1}
\[
H_1(\Ee_2(A),\z)\simeq A/M.
\]
If $|A/\mmm_A|\geq  4$, then by the above argument we have $H_1(\Ee_2(A),\z)=0$. It is 
not difficult to verify that
\[
H_1(\Ee_2(A),\z)\simeq
\begin{cases}
A/\mmm_A^2 &  \text{if $|A/\mmm_A|=2$}  \\
A/\mmm_A &  \text{if $|A/\mmm_A|=3.$}  \\
%0 &  \text{if $|A/\mmm_A|\geq  4$}  \\
\end{cases}
\]
\par (iii) A commutative semilocal ring $A$ is universal for $\GE_2$ if and only if none of the rings 
$\z/2\times\z/2$ and $\z/2 \times \z/3$ is a direct factor of $A/J(A)$ \cite[Theorem 2.14]{menal1979}, 
where $J(A)$ is the Jacobson radical of $A$. Thus if $A$ is semilocal with the above condition, then
by the above corollary
%Theorem~\ref{H1(Y)} and Theorem~\ref{H1} we have
\[
H_1(\Ee_2(A),\z)\simeq A/M.
\]
If $A$ is semilocal such that all its residue fields have more than three elements,
i.e. $|A/\mmm|>3$ for all $\mmm \in {\rm Specm}(A)$, then there is an element $a\in \aa$
such that $a^2-1\in \aa$. Thus $A=M$ and hence
\[
H_1(\Ee_2(A),\z)=0.
\]
(iv) If $A$ is Artinian, then $A\simeq A_{1} \times \cdots \times A_{r}$, where each $A_i$ is an Artinian 
local ring (see \cite[Theorem 8.7]{AM1994}). Hence $\Ee_2(A)\simeq\Ee_2(A_{1})\times\cdots\times\Ee_2(A_{r})$ 
and thus
\[
H_1(\Ee_2(A),\z) \simeq H_1(\Ee_2(A_{1}),\z) \times \cdots \times H_1(\Ee_2(A_{r}),\z).
\]
Now using (ii) we can calculate $H_1(\Ee_2(A),\z)$. 
%For more applications of the above theorem see \cite{B-E2024}.
\end{exa}

\begin{exa}
(i) For a positive square free integer $d$, let $\OO_{-d}$ be the ring of algebraic integers of $\q(\sqrt{-d})$.
%By a well-known result of P. M. Cohn, $\OO_d$ is a $\GE_2$-ring if and only if $d\in\{ -1, -2, -3, -7, -11\}$.
%%\cite[Theorem 6.1]{B-E2024}. 
%Thus $\Ee_2(\OO_d)=\SL_2(\OO_d)$ if and only if $d\in\{ -1, -2, -3, -7, -11\}$.
%He also proved that for $d\neq -1, -2, -3, -7, -11$, we have 
%\[
%H_1(\Ee_2(\OO_{d}),\z)\simeq \z\oplus \z/12.
%\]
Cohn proves that for $-d\neq -1, -2, -3, -7, -11$, $\OO_{-d}$ is discretely normed \cite[\S6]{cohn1966}
and so is universal for $\GE_2$ \cite[Theorem 5.2]{cohn1966}. Since $\OO_{-d}^\times=\{1, -1\}$,
by the above corollary we have
\[
H_1(\Ee_2(\OO_{-d}),\z)\simeq \OO_{-d}/M\simeq \z\oplus \z/12.
\]
%(see \cite[Proposition 4.7 ]{B-E2024}). 
On the other hand Swan has shown that 
\[
H_1(\SL_2(\OO_{-d}),\z)\simeq \begin{cases}
\z\oplus \z \oplus \z/6 \oplus \z/2 & \text{if $d=5$}\\
\z\oplus \z \oplus \z/6 & \text{if $d=6$}\\
\z\oplus \z \oplus \z/12 & \text{if $d=15$}\\
\z & \text{if $d=19$}\\
\end{cases}
\]
(see \cite[Corollaries 11.2, 12.2, 15.2, 16.2]{swan1971}).
Thus the natural map 
\[
H_1(\Ee_2(\OO_{-d}),\z) \arr H_1(\SL_2(\OO_{-d}),\z)
\]
is neither injective nor surjective for $d=5,6$, is injective but not surjective for $d=15$ and
is surjective but not injective for $d=19$.
\par (ii) Let $k$ be a finite field with at least four elements. Then by Example \ref{universal}, 
\[
H_1(\Ee_2(k[X,Y]),\z)=0.
\]
%(seeExample \ref{universal} below), 
But it follows from \cite[Theorem 1.4]{gmv1994} that $H_1(\SL_2(k[X,Y]),\z)$ is an infinite group.
\par (iii) The polynomial ring $A:=\z[X]$ is universal for $\GE_2$ \cite[\S2, Theorem 8.2]{cohn1966}. 
Thus by Theorem \ref{H1(Y)} and Lemma \ref{H0}
\[
H_1(\Ee_2(A), \z)\simeq H_1(\Ee_2(A), Y_\bullet(A^2))\simeq H_1(\SL_2(A), L_\bullet(A^2)).
\]
Observe that by the above corollary 
\[
H_1(\Ee_2(\z[X]), \z)\simeq \z[X]/M\simeq \z[X]/12\z\simeq \z/12 \oplus \bigoplus_{i=1}^\infty \z.
\]
On the other hand, for any positive integer $m$ there is a surjective homomorphism 
\[
\SL_2(\z[X]) \arr F_m, 
\]
$F_m$ a free group with $m$-generators, such that the group of unipotent matrices $U_2(\z[X])$ 
is in its kernel. Since for any ring $A$, $\Ee_2(A) \subseteq U_2(A)$, the natural map 
\[
H_1(\Ee_2(\z[X]), \z) \arr H_1(\SL_2(\z[X]), \z)
\]
has an infinite cokernel. 
\end{exa}

%%%%%%%%%%%%%%%%%%%%%%%%%%%%%%%%%%%%%%%%%%%%%%%%%%%%%%%%%%%%%%%%%%%%%%%%%%
\section{The second homology of \texorpdfstring{$\Ee_2$}{Lg}}\label{s5}
%%%%%%%%%%%%%%%%%%%%%%%%%%%%%%%%%%%%%%%%%%%%%%%%%%%%%%%%%%%%%%%%%%%%%%%%%%

In this section, we study the second homology group of $\Ee_2(A)$, where $A$ is a commutative ring such 
that $2 \in \aa$. The fundamental ideal $I(A)$ plays a central role in this context. More precisely, we 
investigate a surjective map from $H_2(\Ee_2(A), \z)$ onto $I^2(A)$.

We have seen in the proof of Theorem \ref{H1} that the composite
$\bar{I}(A)\arr I(A) \overset{d_{2,0}^2}{\larr} \GG_A\oplus A_\aa$,
takes $\Lan a\Ran\in \bar{I}(A)$, $a\in \aa$, to the element $(\lan a\ran, \overline{3(1-a)})$. 

Let either $2\in \aa$ or $A_\aa=0$. In both cases $\overline{3(1-a)}=0$ in $A_\aa$. So the 
above composite sends $\Lan a\Ran$ to $(\lan a\ran, 0)$. Note that we 
have a surjective map of abelian groups $A/3A \arr A_\aa$. Denote the kernel of 
$d_{2,0}^2: I(A) \arr \GG_A\oplus A_\aa$, by $I^2(A)$, i.e.
\[
I^2(A):=\ker(d_{2,0}^2).
\]
We know that $\bar{I}(A)/\bar{I}^2(A)\simeq \GG_A$ \cite[Theorem 6.1.11]{weibel1994}. 
Consider  the commutative diagram with exact rows 
\[
\begin{tikzcd}
0\ar[r] &\bar{I}^2(A)\ar[r]\ar[d]& \bar{I}(A) \ar[r]\ar[d] &\GG_A\ar[r]\ar[d]&0\\
0 \ar[r] & I^2(A) \ar[r] & I(A) \ar[r] & \GG_A \oplus A_\aa,&
\end{tikzcd}
\]
where the right vertical map is given by $\lan a\ran \mapsto (\lan a\ran, 0)$. By the Snake Lemma
and exact sequence (\ref{I-I}) we obtain the exact sequence
\begin{align}\label{I2-I2}
\bar{I}^2(A) \arr I^2(A) \arr H_1(Y_\bullet(A^2))_{\Ee_2(A)} \arr A_\aa.
\end{align}

\begin{lem}\label{I2--I2}
If $A_\aa=0$, then we have the exact sequence
\[
\bar{I}^2(A) \arr I^2(A) \arr H_1(Y_\bullet(A^2))_{\Ee_2(A)} \arr 0.
\]
\end{lem}
\begin{proof}
This follows immediately from the above exact sequence.
\end{proof}

On the other hand, by an easy analysis of the main spectral sequence we have  the exact sequence
\begin{align}\label{ee}
H_2(\Bb(A),\z) \arr H_2(\Ee_2(A),\z) \arr I^2(A) \arr 0.
\end{align}
We denote the image of $\Lan a\Ran \Lan b\Ran \in \bar{I}^2(A)$ in $I^2(A)$ again by 
$\Lan a\Ran \Lan b\Ran$. Now we lift $\Lan a\Ran \Lan b\Ran\in I^2(A)$ to an element 
of $H_2(\Ee_2(A),\z)$ through the surjective map 
\[
H_2(\Ee_2(A),\z) \arr I^2(A). 
\]

\begin{prp}\label{cycle}
Let $A$ be a commutative ring with $2\in \aa$. If $a, b\in \aa$, then under the surjective map 
\[
H_2(\Ee_2(A),\z) \arr I^2(A),
\]
the homology class of the cycles $([D(a)|D(b)]+R_{ab}-R_a-R_b+R_1)\otimes 1$ maps to 
$\Lan a\Ran \Lan b\Ran\in I^2(A)$, where for any $z\in \aa$,
\begin{align*}
R_z:=&[w|g^{-1}_z]\!-\![w|h^{-1}_z]\!-\![h^{-1}_zD(z)^{-1}|wh^{-1}_z]\!+
\![h_{z^{-1}}|w^{-1}]\!-\![h^{-1}_z|D(z)^{-1}]\!+\![D(z)h_z|D(z)^{-1}]\\
&\!\!\!\!-\![D(z)h_z|h^{-1}_z]\!+\!2[h^{-1}_z|h^{-1}_z]\!+\![h^{-2}_z|h^{-2}_z]\!+
\![D(2)h^{-1}_zD(2)^{-1}|D(2)]\!-\![D(2)|h^{-1}_z]
\end{align*}
with $g_z=\mtxx{0}{1}{-1}{z}$ and $h_z=\mtxx{1}{z^{-1}}{0}{1}$.
\end{prp}
\begin{proof}
We know that $I^2(A) \se I(A) \se \GW(A)=H_0(\Ee_2(A), Z_1(A^2))$. So we need to study the map 
\[
H_2(\Ee_2(A),\z) \arr H_0(\Ee_2(A),Z_1(A^2)).
\]
From the exact sequence $0 \arr Z_1^Y(A^2) \arr Y_1(A^2) \arr  Y_0(A^2) \arr \z \arr 0$
we obtain two short exact sequences
\[
0 \arr Z_1^Y(A^2) \arr Y_1(A^2) \arr  Z_0^Y(A^2) \arr 0, \ \ \ \ \ 
0 \arr Z_0^Y(A^2) \arr  Y_0(A^2) \arr \z \arr 0.
\]
From these we obtain the connecting homomorphisms
\[
\delta_1: H_2(\Ee_2(A),\z) \arr H_1(\Ee_2(A),Z_0^Y(A^2)),
\]
\[
\delta_2: H_1(\Ee_2(A),Z_0^Y(A^2)) \arr H_0(\Ee_2(A),Z_1^Y(A^2)).
\]
Then
\[
\delta_2\circ \delta_1: H_2(\Ee_2(A),\z) \arr H_0(\Ee_2(A),Z_1^Y(A^2))
\]
is the above mentioned map. For $z\in \aa$, the element $\Lan z\Ran\in I(A)$ is represented 
by $[\ ]\otimes\partial_2(X_z-X_1)$, where $X_z:=(\pmb{\infty},\pmb{0},\pmb{z})$. Consider the diagram
\[
\begin{tikzcd}
B_1(\Ee_2(A))\!\otimes_{\Ee_2(A)}\! Z_0^Y(A^2) &\! B_1(\Ee_2(A))\otimes_{\Ee_2(A)}\! Y_1(A^2) 
\!\ar["\id_{B_1}\!\otimes \partial_1"',  l] \ar[d, "d_1\otimes \id_{Y_1}"] & \\
&\!\!\! B_0(\Ee_2(A))\!\otimes_{\Ee_2(A)}\! Y_1(A^2) &\! B_0(\Ee_2(A))\!\otimes_{\Ee_2(A)} 
\! Z_1^Y(A^2).\ar[l, "\id_{B_0}\otimes \inc"']
\end{tikzcd}
\]
Then
\begin{align*}
(\id_{B_0}\otimes \inc)	\bigg([\ ]\otimes\partial_2(X_z-X_1)\bigg)
=&[\ ]\otimes\bigg((\pmb{0},\pmb{z})-(\pmb{\infty},\pmb{z})-(\pmb{0},\pmb{1})+(\pmb{\infty},\pmb{1})\bigg)\\
=&[\ ]\otimes\bigg((g_z-h_z-g_1+h_1)(\pmb{\infty},\pmb{0})\bigg)\\
=&(d_1\!\otimes\!\id_{Y_1(A^2)})\bigg(\!([g^{-1}_z]\!-\![h^{-1}_z]\!-\![g^{-1}_1]\!+
\![h^{-1}_1])\!\otimes\!(\pmb{\infty},\pmb{0})\!\bigg).
\end{align*}
Thus the element 
\begin{align*}
U_z&:=(\id_{B_1}\otimes\partial_1)\bigg(([g^{-1}_z]-[h^{-1}_z]
-[g^{-1}_1]+[h^{-1}_1])\otimes(\pmb{\infty},\pmb{0})\bigg)\\
&=([g^{-1}_z]-[h^{-1}_z]-[g^{-1}_1]+[h^{-1}_1])\otimes 
\partial_1(\pmb{\infty},\pmb{0})
\end{align*}
represents an element of $H_1(\Ee_2(A),Z_0^Y(A^2))$, where maps to $\Lan z\Ran$ through 
$\delta_2$. Now consider the diagram
\[
\begin{tikzcd}
B_2(\Ee_2(A))\otimes_{\Ee_2(A)} \z & B_2(\Ee_2(A))\otimes_{\Ee_2(A)} Y_0(A^2) 
\ar["\id_{B_2}\otimes \varepsilon"',  l] \ar[d, "d_2\otimes \id_{Y_0}"] & \\
& B_1(\Ee_2(A))\otimes_{\Ee_2(A)} Y_0(A^ 2) & B_1(\Ee_2(A))\otimes_{\Ee_2(A)}  
Z_0^Y(A^2).\ar[l, "\id_{B_1}\otimes \inc"']
\end{tikzcd}
\]
Then we have
\begin{align*}
(\id_{B_1}\otimes \inc)(U_z)=&([g^{-1}_z]-[h^{-1}_z]-[g^{-1}_1]
+[h^{-1}_1])\otimes\partial_1(\pmb{\infty},\pmb{0})\\
=&(w-1)([g^{-1}_z]-[h^{-1}_z]-[g^{-1}_1]+[h^{-1}_1])\otimes(\pmb{\infty})\\
=&-[D(z)]+(d_2\otimes\id_{Y_0(A^2)})\bigg((R_z-R_1)\otimes(\pmb{\infty})\bigg).
%\in B_1(\Ee_2(A))\otimes_{\Ee_2(A)} X_0(A^2).
\end{align*}
Note that in our calculation we used the formulas
\[
wg^{-1}_z=D(z)^{-1}h^{-1}_{z^{-1}}wh^{-1}_z, \ \ \ \ \ D(z)h_uD(z)^{-1}=h_{u/z^2},
\]
with $u, z\in\aa$. Now it is easy to see that
\[
(d_2\!\otimes\!\id_{Y_0(A^2)})\bigg(([D(a)|D(b)]\!+\!R_{ab}\!-\!R_a\!-\!R_b\!+\!R_1)
\!\otimes(\pmb{\infty})\bigg)\!=\! (\id_{B_1}\otimes \inc)(U_{ab}\!-\!U_a\!-\!U_b).
\]
Finally, applying $\id_{B_2}\otimes\varepsilon$ we have 
\[
(\id_{B_2}\!\otimes\!\varepsilon)\bigg(\!([D(a)|D(b)]\!+\!R_{ab}\!-\!R_a\!-\!R_b\!+
\!R_1)\!\otimes\!(\pmb{\infty})\!\bigg)
\!=\!([D(a)|D(b)]\!+\!R_{ab}\!-\!R_a\!-\!R_b\!+\!R_1)\!\otimes\! 1.
\]
The element $([D(a)|D(b)]+R_{ab}-R_a-R_b+R_1)\otimes 1$ is a cycle and represents an element 
of $H_2(\Ee_2(A),\z)$ which maps to $\Lan a\Ran \Lan b\Ran\in I^2(A)$. This completes the 
proof of the proposition.
\end{proof}

\begin{rem}
In \cite[\S7]{hutchinson2016}
Hutchinson has constructed a similar cycles as in the above proposition over fields. Let $F$ be a field
with at least four elements. Let $a, b \in F^\times$ and $\lambda \in F^\times$ such that 
$\lambda^2-1\in F^\times$. Then Hutchinson has constructed a cycle
\[
F(a,b)_\lambda\in B_2(\SL_2(F))\otimes_{\SL_2(F)} \z
\]
such that its homology class in $H_2(\SL_2(F),\z)$ is independent of the choice of $\lambda$ and maps to 
$\Lan a \Ran \Lan b \Ran \in I^2(F)$ (see \cite[Proposition 7.1, Lemma 7.5, Theorem 7.8]{hutchinson2016}). 
(Note that the existence of $\lambda \in \aa$ such that $\lambda^2-1 \in A^\times$ implies that $A_\aa=0$.)
\end{rem}

\begin{rem}
Let $A$ be a semilocal ring such that $2\in \aa$ and  all its residue fields have more that four 
elements. Then the group $H_2(\Ee_2(A),\z)$ is generated by the homology class of the cycles 
$([D(a)|D(b)]+R_{ab}-R_a-R_b+R_1)$, $([D(a)|D(b)]-[D(b)|D(a)])$, $a,b \in \aa$, and 
$([E_{12}(x)|E_{12}(y)]-[E_{12}(y)|E_{12}(x)])$, $x,y\in A$:

In fact since $A$ is universal for $\GE_2$ (see Example \ref{universal}), by Proposition \ref{H1(Y)}, 
$H_1(Y_\bullet(A^2))=0$. Thus by Lemma~\ref{I2--I2}, the map $\bar{I}^2(A) \arr I^2(A)$ is surjective. 
This shows that $I^2(A)$ is generated by the elements of the form $\Lan a \Ran\Lan b \Ran$, $a,b\in\aa$. 
Now consider the exact sequence
\[
H_2(\Bb(A), \z) \arr H_2(\Ee_2(A), \z) \overset{\theta}{\arr} I^2(A) \arr 0.
\]
Let $X\in H_2(\Ee_2(A),\z)$ and let $\theta(X)=\sum \Lan a_i \Ran \Lan b_i \Ran$.
If $F(a,b):=([D(a)|D(b)]+R_{ab}-R_a-R_b+R_1)$, then $\theta(X-\sum \overline{F(a_i,b_i)})=0$.
It follows from this that $X-\sum \overline{F(a_i,b_i)}$ sits in the image of $H_2(\Bb(A),\z)$.
Consider the extension $1 \arr \Nn(A) \arr \Bb(A) \arr \Tt(A) \arr 1$, where 
$\Nn(A):=\bigg\{\begin{pmatrix}
1 & b\\
0 & 1
\end{pmatrix}: a\in \aa, b\in A\bigg\}$. 
Since any residue field of $A$ has more than four elements, there is $\lambda \in \aa$ such that 
$\lambda^2-1\in \aa$. This implies that $A_\aa=0$. It follows from this that for any $k\geq 0$, 
$H_k(\aa, A)=0$  \cite[Corollary 3.2]{B-E--2023}. By studying the Lyndon/Hochschild-Serre spectral 
sequence of this extension we obtain the exact sequence
\[
H_2(A,\z)_\aa \arr H_2(\Bb(A), \z) \arr H_2(\Tt(A), \z)\arr 0.
\]
Note that we have the splitting map $\inc_\ast:H_2(\Tt(A), \z) \arr H_2(\Bb(A), \z)$. The image of 
$H_2(A,\z)\simeq A\wedge A$ in $H_2(\Bb(A), \z)$ is generated by the homology classes of the cycles 
\[
G(x,y):=([E_{12}(x)|E_{12}(y)]-[E_{12}(y)|E_{12}(x)])\otimes 1,\ \ \ x,y\in A. 
\]
Moreover, $H_2(\Tt(A), \z)\simeq \Tt(A)\wedge \Tt(A)$ is generated by the homology classes of the cycles 
\[
H(a,b):=([D(a)|D(b)]-[D(b)|D(a)])\otimes 1.
\]
Now the claim follows from the fact that $\ker(\theta)$ is generated by the images of the homology classes 
of the cycles $G(x,y)$ and $H(a,b)$. This proves our claim.
\end{rem}

%%%%%%%%%%%%%%%%%%%%%%%%%%%%%%%%%%%%%%%%%%%%%%%%%%%%%%%%%%%%%%%%%%
\section{The third homology of \texorpdfstring{$\Ee_2$}{Lg}}\label{s6}
%%%%%%%%%%%%%%%%%%%%%%%%%%%%%%%%%%%%%%%%%%%%%%%%%%%%%%%%%%%%%%%%%%

A potential refined Bloch--Wigner exact sequence is an exact sequence that describes the structure of the 
third homology group of $\Ee_2(A)$. Such an exact sequence has been established in a few cases 
\cite{hutchinson-2013}, \cite{hutchinson2017}, \cite{B-E-2023}. We believe that in the present section we 
take a significant step toward a better understanding of such an exact sequence over semilocal rings.

Let $\overline{\RP}(A)$ be the quotient of the free $\z[\GG_A]$-module generated by symbols 
$[a]$, $a\in \WW_A$, by the $\z[\GG_A]$-submodule generated by the elements
\[
[a]-[b]+\lan a\ran\bigg[\frac{b}{a}\bigg]-\lan a^{-1}-1\ran\Bigg[\frac{1-a^{-1}}{1-b^{-1}}\Bigg] +
\lan 1-a\ran\Bigg[\frac{1-a}{1-b}\Bigg],
\]
where $a,b,a/b \in \WW_A$. By a direct computation one can show that the map
\[
\bar{\lambda}_1:\overline{\RP}(A)\arr \II_A^2 \se\z[\GG_A],\ \ \ \ \ \ 
[a]\mapsto\Lan a \Ran\Lan 1-a\Ran,
\]
is a well-defined $\z[\GG_A]$-homomorphism. Let
\[
S_\z^2(\aa):=(\aa \otimes_\z \aa)/\lan a\otimes b+b\otimes a:a,b \in \aa\ran.
\]
If we consider $S_\z^2(\aa)$ as a trivial $\GG_A$-module, then
\[
\bar{\lambda}_2: \overline{\RP}(A) \arr S_\z^2(\aa), \ \ \ \ \ [a] \mapsto a \otimes (1-a),
\]
is a homomorphism of $\z[\GG_A]$-modules. Hutchinson defined the $\z[\GG_A]$-modules 
$\overline{\RP}_1(A)$ and $\overline{\RB}(A)$ as follows
\[
\overline{\RP}_1(A):=\ker(\bar{\lambda}_1:\overline{\RP}(A) \arr \II_A^2),
\]
\[
\overline{\RB}(A):=\ker(\bar{\lambda}_2|_{\overline{\RP}_1(A)}:\RP_1(A) \arr S_\z^2(\aa))
\] 
(see \cite[page 28]{hutchinson-2013}, \cite[Subsection 2.3]{hutchinson2013}).  

From the complex $Y_4(A^2) \arr Y_3(A^2) \arr Z_2^Y(A^2) \arr 0$ we obtain the complex of 
$\GG_A$-modules
\[
Y_4(A^2)_{\Ee_2(A)} \arr Y_3(A^2)_{\Ee_2(A)} \arr Z_2^Y(A^2)_{\Ee_2(A)} \arr 0.
\]
We have seen that 
%the orbits of the action of  $\Ee_2(A)$ on $Y_3(A)$ are represented by
%\[
%\lan a\ran[x]=({\pmb\infty}, {\pmb 0},{\pmb a}, \pmb{ax}),
%\ \ \lan a\ran\in \GG_A, x,\in \WW_A,
%\]and thus 
$Y_3(A^2)_{\Ee_2(A)} $ is a free $\z[\GG_A]$-module generated by the symbols $[x]$, $x\in \WW_A$.
The orbits of the action of  $\Ee_2(A)$ $Y_4(A)$ are represented by
\[  
\lan a\ran[x,y]:= ({\pmb\infty}, {\pmb 0},{\pmb a}, \pmb{ax}, \pmb{ay}),
\ \ \lan a\ran\in \GG_A, x,y,x/y\in \WW_A.
\]
Thus $Y_4(A^2)_{\Ee_2(A)} $ is the free $\z[\GG_A]$-module 
generated by the symbols $[x,y]$, $x,y,x/y\in \WW_A$. It is straightforward to check that 
\[
\overline{\partial_4}([x,y])=[x]-[y]+\lan x\ran\bigg[\frac{y}{x}\bigg]-
\lan x^{-1}-1\ran\Bigg[\frac{1-x^{-1}}{1-y^{-1}}\Bigg]+ 
\lan 1-x\ran\Bigg[\frac{1-x}{1-y}\Bigg].
\]
Thus we obtain a natural map
\[
\eta:\overline{\RP}(A) \arr H_0(\Ee_2(A),Z_2^Y(A^2)).
\]
If $Y_\bullet(A) \arr \z$ is exact in dimension $<4$, then the above map becomes an isomorphism.

Following Coronado and Hutchinson \cite[\S3]{C-H2022} we define the {\it refined scissors 
congruence group} of $A$ as follows:
\[
\RP(A):=H_0(\SL_2(A), Z_2^{L}(A^2)).
\]
It follows from Lemma \ref{ind3}, that
\[
\RP(A)\simeq H_0(\Ee_2(A), Z_2^{Y}(A^2))=Z_2^{Y}(A^2)_{\Ee_2(A)}.
\]
The inclusion $Z_2^Y(A^2) \arr Y_2(A^2)$ induces the natural map 
\[
\lambda_1:=\inc:\RP(A)=Z_2^Y(A^2)_{\Ee_2(A)} \arr Y_2(A^2)_{\Ee_2(A)}\simeq \z[\GG_A].
\]
It is straightforward to check that $\bar{\lambda}_1=\lambda_1\circ\eta$. We denote the kernel 
of $\lambda_1$ by $\RP_1(A)$ and call it the {\it refined scissors congruence group of} $A$. 
Note that we have a natural map 
\[
\overline{\RP}_1(A) \arr \RP_1(A).
\]

Let $Y_\bullet(A)$ be exact in dimension $1$. From the short exact sequence 
\[
0 \arr Z_2^Y(A^2)\arr Y_2(A^2) \overset{\partial_2^Y}{\arr} Z_1^Y(A^2) \arr 0,
\]
we obtain the (long) exact sequence
\[
H_1(\Ee_2(A), Y_2(A^2)) \overset{\bar{\partial}_2^Y}{\arr} H_1(\Ee_2(A), Z_1^Y(A^2)) 
\arr \RP(A) \overset{\lambda_1}{\larr} \z[\GG_A].
\]
Since
\[
H_1(\Ee_2(A), Y_2(A^2)) \simeq H_1\bigg(\Ee_2(A), \bigoplus_{\lan a \ran \in \GG_A} 
\Ind_{\mu_2(A)}^{\Ee_2(A)}\z\bigg)\simeq \mu_2(A) \otimes_\z \z[\GG_A],
\]
we have the exact sequence
\[
\mu_2(A) \otimes_\z \z[\GG_A] \overset{\bar{\partial}_2}{\arr} \mathrm{E}_{2,1}^1 \arr 
\RP_1(A) \arr 0.
\]
From the commutative diagram
\[
\begin{tikzcd}
Y_2(A^2) \ar[r, "\partial_2^Y"]\ar[d, "\partial_2"] & Z_1^Y(A^2) \ar[d, "\inc"]\\
Y_1 (A^2) \ar[r, equal] & Y_1 (A^2)
\end{tikzcd}
\]
we obtain the commutative diagram
\[
\begin{tikzcd}
\mu_2(A)\otimes_\z\z[\GG_A]\ar[r, "\bar{\partial}_2"] \ar[d] & \mathrm{E}_{2,1}^1\ar[r] 
\ar[d, "d_{2,1}^1"] &\RP_1(A)\larr 0\\
\mu_2(A) \ar[r, equal] & \mu_2(A).
\end{tikzcd}
\]
It is not difficult to see that $(d_{2,1}^1\circ \bar{\partial}_2^Y)(b\otimes\lan a\ran)=b$. Thus
from the above commutative diagram we obtain the exact sequence
\[
\mu_2(A) \otimes_\z \II_A \arr \mathrm{E}_{2,1}^2 \arr \RP_1(A) \arr 0.
\]
Consider the differential $d_{2,1}^2:E_{2,1}^2 \arr H_2(\Bb(A),\z)$. In the rest of this section
we will assume that the surjective map $p:\Bb(A) \arr \Tt(A)$ induces the isomorphism
\[
H_2(\Bb(A), \z) \overset{p_\ast}{\simeq} H_2(\Tt(A),\z).
\]
For example over semilocal rings we have the following result.

\begin{prp}\label{iso-hut}
Let $A$ be a semilocal ring. For any maximal ideal $\mmm\in {\rm Specm}(A)$, let
either $A/\mmm$ be infinite or if $|A/\mmm|=p^{d}$, then $(p-1)d>2(n+1)$ ($(p-1)d>2n$  
when $A$ is a domain). Then $H_k({\rm T}(A), \z)\simeq H_k({\rm B}(A),\z)$ for all $k\leq n$.
\end{prp}
\begin{proof}
If $A$ is a domain, then the proof is very similar to the proof of a similar result over local domains 
presented in \cite[\S3]{hutchinson2017}. Just there one should replace the maximal ideal of the local 
domain (denoted by $\mathcal{M}$ in \cite[\S3]{hutchinson2017}) with the Jacobson radical $J(A)$ of $A$. 

Now let $A$ be a general semilocal ring with the given condition. Then similar to the proof of 
\cite[Proposition 3.19]{hutchinson2017}, we can show that $H_n(\Tt(A), k)\simeq H_n(\Bb(A), k)$, where $k$ 
is a prime field and $(p -1)d > 2n$. Now the claim follows from \cite[Lemma 2.3]{mirzaii2017}.
\end{proof}

Thus let $H_2(\Bb(A), \z) \overset{p_\ast}{\simeq} H_2(\Tt(A),\z)$. Then the composite
\[
\mu_2(A)\otimes_\z\II_A \arr \mathrm{E}_{2,1}^2\overset{d_{2,1}^2}{\larr} H_2(\Bb(A),\z)\simeq 
H_2(\Tt(A),\z)\simeq\aa\wedge \aa
\]
is given by $b \otimes \Lan a\Ran \mapsto b\wedge a$ (see \cite[Lemma 4.1]{B-E--2023}). Thus we have 
the commutative diagram
\[
\begin{tikzcd}
&\mu_2(A) \otimes_\z \II_A\ar[r] \ar[d] & \mathrm{E}_{2,1}^2 \ar[r] \ar[d,"d_{2,1}^2"] 
&\RP_1(A)\ar[r]\ar[d] & 0\\
0 \ar[r] &\mu_2(A)\wedge \aa \ar[r] & \aa \wedge \aa \ar[r] & \displaystyle 
\frac{\aa\wedge \aa}{\mu_2(A)\wedge \aa} \ar[r] & 0.
\end{tikzcd}
\]
Consider the injective map
\[
\alpha: \frac{\aa\wedge \aa}{\aa\wedge \mu_2(A)} \arr S_\z^2(\aa), \ \ \ \ 
a\wedge b \mapsto 2(a\otimes b).
\]
One can show that the composite map 
\[
\overline{\RP}_1(A) \arr \RP_1(A) \arr \frac{\aa\wedge \aa}{\mu_2(A)\wedge \aa} \arr S_\z^2(\aa)
\]
coincides with $\bar{\lambda}_2|_{\overline{\RP}_1(A)}: \overline{\RP}_1(A) \arr S_\z^2(\aa)$. 
We denote the kernel of 
\[
\RP_1(A) \arr \frac{\aa\wedge \aa}{\mu_2(A)\wedge \aa} \arr S_\z^2(\aa)
\]
with $\RB(A)$ and call it the {\it refined Bloch group of} $A$. Note that we have a natural map
\[
\overline{\RB}(A) \arr \RB(A).
\]
Now from the above commutative diagram we obtain a surjective map $\mathrm{E}_{2,1}^\infty \two \RB(A)$.
By an easy analysis of the main spectral sequence we obtain the surjective map 
\[
H_3(\Ee_2(A),\z) \arr \mathrm{E}_{2,1}^\infty.
\]
Combining these two maps we obtain a surjective map 
\[
H_3(\Ee_2(A),\z) \arr \RB(A).
\]
This map factors through $H_3(\PEe_2(A),\z)$, where 
\[
\PEe_2(A):=\Ee_2(A)/\{\pm I_2\}.
\]
We have proved the following version of a refined Bloch-Wigner exact sequence in \cite{B-E-2024}.

\begin{thm}\label{Proj-BW}
Let $A$ be a semilocal ring such that there is a ring homomorphism $A\arr F$, $F$ a field, where 
$\mu(A) \arr \mu(F)$ is injective. For any maximal ideal $\mmm \in {\rm Specm}(A)$ let either 
$A/\mmm$ be infinite or if $|A/\mmm|=p^{d}$, then $(p-1)d>8$ ($(p-1)d>6$ when $A$ is a domain). 
If either $-1\in \aa^2$ or $|\GG_A|\leq 4$, then we have the exact sequence
\[
0\arr \tors(\widetilde{\mu}(A),\widetilde{\mu}(A)) \arr H_3(\PEe_2(A),\z) \arr \RB(A) \arr 0,
\]
where $\widetilde{\mu}(A)=\mu(A)/\mu_2(A)$.
\end{thm}
\begin{proof}
See \cite[Theorem 3.1, Theorem 3.2, Corollary 3.3]{B-E-2024} and use Proposition \ref{iso-hut} 
and \cite[Lemma 2.8]{B-E-2024}.
\end{proof}

Let $\Aa$ be a finite cyclic group. If $2\mid |\Aa|$, let $\Aa^\sim$ denote the unique non-trivial 
extension of $\Aa$ by $\z/2$. If $2\nmid|\Aa|$, we define $\Aa^\sim:=\Aa$. Thus if $n=|\Aa|$, then
\[
\Aa^\sim \simeq \z/\gcd(n^2, 2n).
\]

Let $\Aa$ be an ind-cyclic group, i.e. $\Aa$ is direct limit of its finite cyclic subgroups. Then 
by passing to the limit we can define $\Aa^\sim$. Note that if $\Aa$ has no $2$-torsion, then 
$\Aa^\sim=\Aa$. There is always a canonical injective homomorphism $\Aa \harr \Aa^\sim$ whose 
composition with the projection $\Aa^\sim \arr \Aa$ coincides with multiplication by $2$. 

The connection between the third homology groups of $\Ee_2(A)$ and $\PEe_2(A)$ is given in 
\cite[Proposition 5.1]{B-E---2024}.

\begin{thm}\label{SL-PSL--1}
Let $A$ be a semilocal ring such that there is a ring homomorphism $A\arr F$, $F$ a field, where 
$\mu_2(A) \simeq \mu_2(F)$. If all the residue fields of $A$ have more than three elements, 
then we have the exact sequence 
\[
0 \arr \mu_2(A)^\sim \arr\frac{H_3(\Ee_2(A),\z)}{\mu_2(A)\otimes_\z H_2(\Ee_2(A),\z)} \arr 
H_3(\PEe_2(A),\z) \arr 0,
\]
where the map $\mu_2(A)\otimes_\z H_2(\Ee_2(A),\z) \arr H_3(\Ee_2(A),\z)$ is induced by the
product map $\mu_2(A) \times \Ee_2(A) \arr \Ee_2(A)$, $(-1, X) \mapsto -X$.
\end{thm}
\begin{proof}
We may assume that $\char(F)\neq 2$. Since for any maximal ideal $\mmm\in {\rm Specm}(A)$, 
$|A/\mmm|\geq 4$, there is $a\in \aa$ such that $a^2-1\in \aa$. It follows from this that 
$M=A$. Thus by Theorem \ref{H1} we have $H_1(\Ee_2(A),\z)=0$. Now by \cite[Theorem 3.1]{B-E---2024} 
we obtain the desired exact sequence.
\end{proof}

The above two theorems are results of our attempt to find an answer to the following question 
of Coronado and Hutchinson \cite[page 3]{C-H2022}.

\begin{question}\label{Q2}
Is it true that for any local domain $A$ one has the exact sequence
\[
0\arr \tors(\mu(A),\mu(A))^\sim \arr H_3(\SL_2(A),\z)\arr\RB(A)\arr 0?
\]
\end{question}

The main result of this section is as follows, which generalizes the main result of 
\cite[Theorem 5.1, Corollary 5.2]{B-E-2023}.

\begin{thm}\label{HBW}
Let $A$ be a semilocal ring such that there is a ring homomorphism $A\arr F$, $F$ a field, where 
$\mu(A) \arr \mu(F)$ is injective. For any maximal ideal $\mmm \in {\rm Specm}(A)$ let either 
$A/\mmm$ be infinite or if $|A/\mmm|=p^{d}$, then $(p-1)d>8$ ($(p-1)d>6$ when $A$ is a domain).
If either $-1\in \aa^2$ or $|\GG_A|\leq 4$, then we have the exact sequence
\[
0\arr \tors(\mu(A), \mu(A))^\sim \arr \frac{H_3(\Ee_2(A),\z)}{\mu_2(A)\otimes_\z I^2(A)}
\arr\RB(A) \arr0,
\]
where $(-1)\otimes \Lan a\Ran\Lan b\Ran \in \mu_2(A)\otimes_\z I^2(A)$ maps in $H_3(\Ee_2(A),\z)$ 
to the class of the cycle which is obtained from the shuffle product of $[-1]\in B_1(\mu_2(A))_{\mu_2(A)}$ 
with the cycle $[D(a)|D(b)]+R_{ab}-R_a-R_b+R_1 \in B_2(\Ee_2(A))_{\Ee_2(A)}$.
\end{thm}
\begin{proof}
Consider the surjective maps 
\[
\frac{H_3(\Ee_2(A),\z)}{\mu_2(A)\otimes_\z H_2(\Ee_2(A),\z)} \two H_3(\PEe_2(A),\z) \two \RB(A).
\]
Let $\KK$ be the kernel of this composite. Then we have the commutative diagram with exact rows
\[
\begin{tikzcd}
0\ar[r] &\KK\ar[r]\ar[d]& \displaystyle\frac{H_3(\Ee_2(A),\z)}{\mu_2(A)\otimes_\z H_2(\Ee_2(A),\z)} 
\ar[r]\ar[d, two heads] &\mathcal{\RB}(A)\ar[r]\ar[d, equal]&0\\
0 \ar[r] &\tors(\widetilde{\mu}(A),\widetilde{\mu}(A))\ar[r] &H_3(\PEe_2(A),\z)\ar[r] &\RB(A)\ar[r]& 0
\end{tikzcd}
\]
(see Theorem \ref{Proj-BW}). From this diagram and Theorem \ref{SL-PSL--1} we obtain the exact sequence
\[
0 \arr \mu_2(A)^\sim \arr \KK \arr \tors(\widetilde{\mu}(A), \widetilde{\mu}(A)) \arr 0.
\]
This clearly implies that $\KK\simeq \tors(\mu(A), \mu(A))^\sim$ and thus we have the exact sequence
\[
0 \arr\tors(\mu(A), \mu(A))^\sim \arr\frac{H_3(\Ee_2(A),\z)}{\mu_2(A)\otimes_\z H_2(\Ee_2(A),\z)} 
\arr \RB(A) \arr 0.
\]
By Proposition \ref{iso-hut}, $H_2({\rm B}(A),\z)\simeq H_2({\rm T}(A),\z)$. Under the composition 
\[
{\rm T}(A)\wedge {\rm T}(A)\overset{\simeq}{\larr} H_2({\rm T}(A),\z) \overset{\simeq}{\larr}  
H_2({\rm B}(A),\z)\arr H_2(\Ee_2(A),\z)
\]
$D(a)\wedge D(b)$ maps to the homology class of the cycle $([D(a)|D(b)]-[D(b)|D(a)])$, $a, b\in\aa$.
Now from the commutative diagram
\[
\begin{tikzcd}
\mu_2(A)\otimes_\z ({\rm T}(A)\wedge {\rm T}(A))\ar[r] \ar[d] & \mu_2(A)\otimes_\z H_2({\rm T}(A),\z) 
\ar[r]\ar[d] & \mu_2(A)\otimes_\z H_2(\Ee_2(A),\z) \ar[d]\\
{\rm T}(A)\wedge {\rm T}(A)) \wedge {\rm T}(A) \ar[r] & H_3({\rm T}(A),\z) \ar[r] & H_3(\Ee_2(A),\z)
\end{tikzcd}
\]
it follows that the right vertical map is given by
\[
(-1)\otimes (D(a)\wedge D(b))\mapsto D(-1)\wedge D(a)\wedge D(b).
\]
Note that here we denoted the image of $D(-1)\wedge D(a)\wedge D(b)\in \bigwedge_\z^3 \Tt(A)$ in 
$H_3(\Ee_2(A),\z)$ again by $D(-1)\wedge D(a)\wedge D(b)$. If $\ii^2=-1$, then 
\[
D(-1) \wedge D(a)\wedge D(b)=2(D(\ii)\wedge D(a)\wedge D(b))=0.
\]
Let $|\GG_A|\leq 4$. If $-1 \in \aa^2$, we proceed as above. If $-1 \notin \aa^2$, then
$\GG_A=\{1, -1, u, -u\}$ for some $u\in \aa$. So we may take $a, b \in \{1, -1, u, -u\}$ and
in any case
\[
D(-1)\wedge D(a)\wedge D(b)=0.
\]
Thus
%\[
$\displaystyle\frac{H_3(\Ee_2(A),\z)}{\mu_2(A)\otimes_\z H_2(\Ee_2(A),\z)} 
=\frac{H_3(\Ee_2(A),\z)}{\mu_2(A)\otimes_\z I^2(A)}$.
%\]
This completes the proof of the theorem.
\end{proof} 

In order to have a positive answer to Question \ref{Q2}, we should have a positive answer to
the following question, which was raised in \cite{B-E---2024}.

\begin{question}\label{Q1}
Let $A$ be a local domain such that $|A/\mmm_A|\neq 2$. Is the product map
\[
\rho_\ast:\mu_2(A)\otimes_\z H_2(\Ee_2(A),\z) \arr H_3(\Ee_2(A),\z)
\]
trivial?
\end{question}

We do not know the answer of the above question even over a general infinite field. For more on 
this we refer the interested readers to \cite[Example 4.4, Remark 5.2, Proposition 5.4]{B-E---2024}.

\begin{rem}
Let $A$ be a local domain such that $|A/\mmm_A| \neq 2, 3,4,5,7,8,9,16,27,32,64$. Then the results 
of this paper suggest that to  give a positive answer to Question \ref{Q2}, one may show that
\par (i) we have the projective Bloch-Wigner exact sequence
\[
0 \arr \tors(\widetilde{\mu}(A),\widetilde{\mu}(A)) \arr H_3(\PEe_2(A),\z) \arr \RB(A) \arr 0,
\]
\par (ii) the product map $\rho_\ast:\mu_2(A)\otimes_\z H_2(\Ee_2(A),\z) \arr H_3(\Ee_2(A),\z)$ 
is trivial.
\end{rem}

\begin{rem}
Although significant progress has been made toward understanding the low-dimensional 
homology of $\Ee_2(A)$ (see the list of references below), many problems remain open for 
future research. 

For the first homology, we hope to obtain a deeper understanding of the kernel of the 
surjective map $A/M \arr H_1(\Ee_2(A), \z)$ (from Theorem~\ref{H1}) for an arbitrary 
commutative ring. This requires a more detailed analysis of the structure of the 
group $H_1(Y_\bullet)$.

Although we have a reasonable understanding of the second homology of $\Ee_2(A)$ over 
semilocal rings, determining the precise structure of this group over general rings 
remains an open problem. For results concerning the structure of the group 
$H_2(\Ee_2(A),\z)$ for $A=\z[1/n]$, see \cite{hutchinson2016} and \cite{MRV2025} and
for $A$ a discrete valuation rings see \cite{HMM2022}.

It appears that the third homology of $\Ee_2(A)$ over semilocal rings is best described by 
a refined Bloch--Wigner exact sequence. As explained in the previous remark, a positive 
answer to Problem~\ref{Q2} depends on a positive answer to Problem~\ref{Q1} and on the 
existence of a projective refined Bloch--Wigner exact sequence. We believe that both 
questions should have affirmative answers, provided that the residue fields of the local 
rings are sufficiently large. It is worth noting that even over a general infinite field, 
a definitive result is still lacking.
\end{rem}

%%%%%%%%%%%%%%%%%%%%%%%%%%%%%%%%%%%%%%%%%%%%%%%%%%%%%%%%

\end{document}